# Joint probability density functions of random trajectories through a box

*Gregory T. Clement*

Department of Biomedical Engineering, Cleveland Clinic Foundation, 9500 Euclid Ave/ND 20, Cleveland, Ohio 44195: clemeng@ccf.org

## Abstract

A wide range of physical problems can be described by randomly-oriented linear trajectories, including any system of objects, organisms, particles, or rays that follow a linear path. Dependent upon the particular random variables that define translation and direction, a description of the probabilistic behavior through a static volume can reveal evidence on an object's physical properties. Yet more information may be available in cases where this behavior is a function of volume shape. In this study two types of random trajectories are considered as they pass through a box of arbitrary relative dimension. One type defines trajectories from uniform random selection of a spatial location paired with a directional vector; the other is formed from a uniformly distributed position vector on a surface. The joint probability distributions for trajectory length as a function of box position are formulated and then examined for different size boxes, and their physical representations discussed.

## 1. Introduction

A wide range of objects, organisms, particles, and rays follow a linear path whose orientation is randomly distributed. The probabilistic behavior of such trajectories passing through a confined volume can reveal unique - and seemingly underutilized – geometrically-dependent values that can give insight on the particle type or, conversely, the geometry being trajected. For example, the mean distance across a volume forms a set of fundamental constants for a given geometry [1], which has led to such applications as the probabilistic determination partial size distributions [2], [3].

Interpretation also depends upon the particular type of randomness that governs a trajectory, and indeed a valid representation of a specific physical behavior requires that the criteria for randomness be accurately described. For the case of a cubic volume, several constructs of random behavior were first considered over 40 years ago by Coleman [4], who later solved similar cases for the probability density functions (PDF) through a box of arbitrary dimension [5]. The present work likewise investigates PDFs and expected distances of randomly-oriented vectors through a box, concentrating on two specific cases that may have application in wave and particle detection along a surface. While Colman [5] elegantly provides an exact form solution to the distribution of path lengths, the current work investigates joint distributions formulated as a function of distance and exit position. Interest in these distributions stems from their potential relevance to experiments involving a finite number of detectors placed at fixed locations on a volume surface. Of course the volume need not be physical, but rather may define positions that exploit the probabilistic behavior of a given geometry.

Two types of random trajectories are considered; the first being potentially relevant to particle or ray propagation and the second which may apply to cases of point-like wave sources. **Case I** defines trajectories from uniform random selection of a spatial location paired with a directional vector, also randomly and uniformly chosen. Such trajectories comprise the most unconstrained set of linear paths that might happen to traverse a given volume. **Case II** concerns a line whose entrance into the box is randomly and uniformly distributed, but whose probability of exiting any given point on the surface is unity. Such behavior, for example, might describe a point-like wave source appearing randomly on the surface but who's radiated field exits over the entire surface (e.g. the emission from a collapsing bubble that forms on an interface).

Although the following formulation has been generalized, much impetus for this work stems from our ongoing studies on the flow of micro-bubbles through blood vessels [6], [7], where analysis of 3-dimensional ultrasound signals as a function of time and bubble density within a given volume reveals dependency on volume shape that we hypothesize can be ascribed to probabalistic behavior. Interestingly, the problem has also recently found potential application in our studies on the trajectories of ultrasound wavevectors, as applied toward transmission ultrasound imaging of the human brain [8].

# 2. Formulation

## 2.1 Entrance

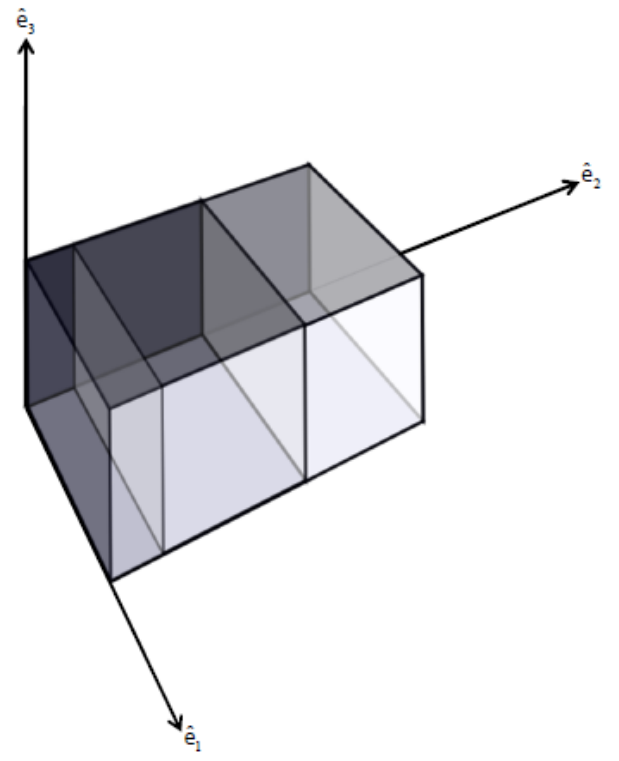

Figure 1. Box configuration.

The goal of the following section is to describe the joint distribution for length and exiting surface position of trajectories through a box. Cases of randomness are considered include: Case I, where trajectories are formed by uniform random selection of a point in space paired with a random directional vector, and Case II, where trajectories are defined by a single random surface point extending radially outward to all surface points throughout the box. Trajectories defined from two uniformly and randomly distributed points the box surface will also be used briefly for comparison in Section 3. However, as this case has been well described, at least for the case of a cube [9], [10] a detailed development is omitted here.

All cases noted above may be described in terms of trajectories that happen to intersect a square box (cuboid), of length $X_1$, depth $X_2$ and height $X_3$, arbitrarily placed with three edges coinciding with the Cartesian axes depicted in Figure 1.

For (I), randomly oriented lines are defined from a directional vector

$$\mathbf{n} = \mathrm{n}\left(\frac{\tilde{x}_1}{\mathrm{r}}\hat{\mathbf{e}}_1 + \frac{\tilde{x}_2}{\mathrm{r}}\hat{\mathbf{e}}_2 + \frac{\tilde{x}_3}{\mathrm{r}}\hat{\mathbf{e}}_3\right);\ \mathrm{r}=\sqrt{\tilde{x}_1{}^2+\tilde{x}_2{}^2+\tilde{x}_3{}^2},\ \mathrm{r} \le 1;\ -1 \le \tilde{\mathrm{x}}_\mathrm{i} \le 1,\ \mathrm{i}=1,2,3. \tag{1}$$

paired with a randomly-selected and point of intersection with the Cartesian plane over $x_i$ =0. The criterion that $\mathrm{r} \le 1$ assures a uniform distribution [11]. If such intersections are uniformly distributed over the entire plane, a vector that traverses the box will enter at uniformly distributed locations over its surface. Hence, the probability of entering the box on a given side for case (I), (II), and (III) are identical.

Despite the arbitrary relative lengths of ($X_1$, $X_2$, $X_3$), it may be observed that substantial symmetry exists in the numeric description of the problem. Vectors entering the box through a given side will exit either through (a) the opposing face, or (b) one of four adjacent sides. Introducing the indices ( i, j, k) to represent any *even* permutation of (1,2,3) with $\mathrm{i} \ne \mathrm{j} \ne \mathrm{k}$ will allow reduction of the discussion to cases of exit through opposing or adjacent faces. Continuing under this notation, the probability of entering a given side can be written as the ratio of the side area to the total surface area:

$$\mathrm{P}_{\mathrm{X_iX_j}} = \mathrm{X_iX_j}/(2\mathrm{X_iX_j}+2\mathrm{X_iX_k}+2\mathrm{X_jX_k}). \tag{2}$$

## 2.2 Case I: Random orientation

### *2.2.1. Opposing faces: Distribution of exit location*

For a trajectory entering the box through $\mathrm{x_j} = 0,\ \mathbf{r} = \mathrm{x_i}\hat{\mathbf{e}}_\mathrm{i} + \mathrm{x_k}\hat{\mathbf{e}}_\mathrm{k}$, and exiting on its opposing face $\mathrm{x'_j} = \mathrm{X_j},\ \mathbf{r}' = \mathrm{x}'_\mathrm{i}\hat{\mathbf{e}}_\mathrm{i} + \mathrm{X_j}\hat{\mathbf{e}}_\mathrm{i} + \mathrm{x}'_\mathrm{k}\hat{\mathbf{e}}_\mathrm{k}$,, the PDF of the distance across the box can be solved by first considering the distribution over the entire plane $\left(x_i', X_j, x_k'\right)$. The vector from the entry point to this plane is given by

$$\mathbf{n} = (x_i' - x_i)\hat{\mathbf{e}}_i + X_j\hat{\mathbf{e}}_j + (x_k' - x_k)\hat{\mathbf{e}}_k \cdot \tag{3}$$

Equating the j-dimension of (3) with (1), the vector length is written

$$n = X_j \frac{r}{\tilde{x}_j}. \tag{4}$$

Substituting (4) into the i-dimension of (1) and using (3),

$$x_i' = x_i + \tilde{x}_i \frac{X_j}{\tilde{x}_j} \cdot \tag{5}$$

In this form, $x'_i$ is described in terms of three independent variables having density distributions

$$\left.\begin{array}{l} f_{x_i}(x_i) = \begin{cases} \dfrac{1}{X_i}; & 0 \le x_i \le X_i \\ 0; & \textit{Otherwise} \end{cases} \\ f_{\tilde{x}}(\tilde{x}_i) = \begin{cases} \dfrac{1}{2}; & -1 \le \tilde{x}_i \le 1 \\ 0; & \textit{Otherwise} \end{cases} \end{array}\right\} i=1,2,3 \tag{6}$$

The PDF of (5) is determined by combining (6) according to relations (A.1a) and (A.2) in Appendix A:

$$f_{x'_i}(x'_i) = \frac{1}{X_j} f_{x_i}(x_i) \otimes \int f_{\tilde{x}}\left(\frac{x_i w}{X_j}\right) f_{\tilde{x}}(w)\, w\, d\, w, \tag{7}$$

which provides the distribution of $x_i'$ along $x_j' = X_j$.

The joint density $f\left(x_i', x_k'\right)$ can now be solved by first treating $x_i'$ as a constant over $0 \le x_i' \le 1$ and solving for $f\left(x_i', x_k'\right)$. The result can then be combined with (6) to form $f\left(x_i', x_k'\right) = f\left(x_i'|x_k'\right) f_{x_i}\left(x_i'\right)$. Rearranging (5), the PDF of the term $\alpha_j = (X_j/\tilde{x}_j)|x_i' = (x_i' - x_i)/\tilde{x}_i$ is solved using (A.2)

$$f_{\alpha_j}(\alpha_j) = \int f_{x_i}(x_i' - \alpha_j w) f_{\tilde{x}}(w)\, w dw \cdot \tag{8}$$

Solving (1) and (3) for $x_k'$ and substituting $\alpha_j$,

$$x_k' = x_k + \tilde{x}_k \alpha_j, \tag{9}$$

so that by (A.1a), (A.3),and (8)

$$f_{x'_k}(x'_k|x'_i) = f_{x_k}(x_k) \otimes \int \frac{f_{\tilde{x}_k}(x_k/w) f_\alpha(w)}{w} dw \cdot \qquad (10)$$

Equations (7) and (10) are then combined to give:

$$f_{jj}(x'_i, x'_k)\,|_{x_j=0}^{x'_j=X_j} = f_{x'_k}(x'_k \mid x'_i)\, f_{x'_i}(x'_i)\,|_{x_j=0}^{x'_j=X_j}, \qquad (11)$$

the PDF for exiting at $(x'_i, x'_k)$ on the side $x'_j = X_j$ when a trajectory enters on the side $x_j = 0$. By symmetry it must also be the case that $f_{jj}(x'_i, x'_k)\,|_{x_j=0}^{x'_j=X_j} = f_{jj}(x'_i, x'_k)\,|_{x_j=X_j}^{x'_j=0}$, such that even permutations describe three unique PDFs representing entrance through any of the six sides.

*2.2.2. Opposing faces: Joint distribution of exit location and path length*

The joint distribution is solved through a change of variables between the five independent variables $(x_i, x_k, \tilde{x}_i, \tilde{x}_j, \tilde{x}_k)$ and $(n, r, x'_i, x'_k, w_i)$; $w_i = x_i$,

$$f(n, r, x'_i, x'_k, w_i) = \frac{f_{x_i}(x_i) f_{x_k}(x_k) f_{\tilde{x}}(\tilde{x}_i) f_{\tilde{x}}(\tilde{x}_j) f_{\tilde{x}}(\tilde{x}_k)}{|J|} \qquad (12)$$

where components of the Jacobian, J, are obtained from the first partial derivatives of (4) and (5) along with $\partial w_i / \partial x_i = 1$. Appendix B supplies details on the solution for J, whose solution is given by (B.1),

$$J = \frac{n^2 \sqrt{n^2 - (x'_i - w_i)^2 - X_j^2}}{X_j r^2}. \qquad (13)$$

Changing variables in (12) and integrating over r and $w_i$, leads to the desired distribution,

$$f_{njj}(n, x'_i, x'_k) = \qquad (14)$$

$$X_j \int_{r, w_i} \frac{f_{x_i}(w_i) f_{x_k}(x'_k - \sqrt{n^2 - (x'_i - w_i)^2 - X_j^2})\, f_{\tilde{x}_i}(\frac{r(x'_i - w_i)}{n})\, f_{\tilde{x}_j}(X_j \frac{r}{n})\, f_{\tilde{x}_k}(r\sqrt{1 - \frac{(x'_i - w_i)^2 - X_j^2}{n^2}})}{n^2 \sqrt{n^2 - (x'_i - w_i)^2 - X_j^2}} r^2 dr dw_i$$

where the symmetry described for (11) applies so that the even permutations generalize to all six sides.

*2.2.3. Adjacent faces: Distribution of exit location*

A line entering the box through the side along $x_j = 0$ is now considered for a trajectory exiting on an adjacent side $x'_k = 0$. Using $n / r = -x_k / \tilde{x}_k$, the $x'_i$ is described by

$$x_i' = x_i - \tilde{x}_i \frac{x_k}{\tilde{x}_k} \cdot \tag{15}$$

Letting $\beta_i = \frac{\tilde{x}_i}{\tilde{x}_k}$, by (A.2),

$$f_{\beta_i}(\beta_i) = \int f_{\tilde{x}}(\beta_k w) f_{\tilde{x}}(w)\, w\, d\, w \, , \tag{16}$$

then by (A.1b) and (A.3),

$$f_{x_k'}(x_k' | x_i') = f_{x_k}(x_k) \otimes \int \frac{f_{\tilde{x}_k}(x_k / w) f_\beta(w)}{w} dw$$

Mirroring arguments leading to (11), the PDF for the position $(x_i', x_j')\,|_{x_j=0}^{x_k'=0}$ is first determined. Treating $x_j'$ as a constant and writing $\alpha_i = n / r = (x_i' - x_i) / \tilde{x}_i$, the PDF of $x_j' = x_j \left( x_i' - x_i \right) / \tilde{x}_i$ can be found by utilizing (8), and by letting $x_j' = x_j \alpha_i$. Then by A.3

$$f_{x_j'}(x_j' | x_i') = \int \frac{f_{\alpha_i}(x_j' / s) f_{\tilde{x}_j}(s)}{s} ds \cdot \tag{17}$$

By (17) and (12),

$$f(x_j', x_i') = f_{x_j'}(x_j' | x_i') f_{x_i'}(x_i') \cdot \tag{18}$$

Once again considering symmetry, (18) may be equated to its opposing adjacent side $f_n(x_i', x_j')\,|_{x_j=0}^{x_k'=0}$ $= f(x_i', x_j')\,|_{x_j=0}^{x_k'=X_k}$ Moreover, by symmetry of the entering face with its opposite side,

$$\begin{aligned} f(x_i', x_j')\Big|_{x_j=0}^{x_k'=0} &= f(x_i', x_j')\,|_{x_j=0}^{x_k'=X_k} \\ &= f(x_i', x_j')\,|_{x_j=X_j}^{x_k'=0} \; ; \quad \textit{all permutations} \\ &= f(x_i', x_j')\,|_{x_j=X_j}^{x_k'=X_k} \end{aligned} \tag{19}$$

Taken over *all* permutations of (i,j,k) gives six unique PDFs. Equations (18) along with (11) form the set of 9 PDFs for the problem.

*2.2.4. Adjacent faces: Joint distribution of exit location and path length*

For adjacent faces, the joint distribution will be solved through a change of variables between $(x_i, x_k, \tilde{x}_i, \tilde{x}_j, \tilde{x}_k)$ and $(n, x_i', x_j', \zeta_k, w_i)$; $w_i = x_i$, $\zeta_k = \frac{x_k}{\tilde{x}_k}$,

$$f(\mathrm{n}, \mathrm{x}_\mathrm{i}', \mathrm{x}_\mathrm{j}', \zeta_\mathrm{k}, \mathrm{w}_\mathrm{i}) = \frac{f_{\mathrm{x}_\mathrm{i}}(\mathrm{x}_\mathrm{i}) f_{\mathrm{x}_\mathrm{k}}(\mathrm{x}_\mathrm{k}) f_{\tilde{\mathrm{x}}_\mathrm{i}}(\tilde{\mathrm{x}}_\mathrm{i}) f_{\tilde{\mathrm{x}}_\mathrm{j}}(\tilde{\mathrm{x}}_\mathrm{j}) f_{\tilde{\mathrm{x}}_\mathrm{k}}(\tilde{\mathrm{x}}_\mathrm{k})}{|\mathrm{J}|}. \tag{20}$$

Components of the Jacobian, J, are again obtained from the first partial derivatives of (4) and (5) along with $\partial \mathrm{w}_\mathrm{i} / \partial \mathrm{x}_\mathrm{i} = 1, \partial \zeta_z / \partial \mathrm{x}_z = 1 / \tilde{\mathrm{x}}_z$ and $\partial \zeta_z / \partial \mathrm{x}_z = -\mathrm{x}_z / \tilde{\mathrm{x}}_z^{\,2}$. The solution, provided in (B.2), is

$$\mathrm{J} = \frac{\zeta_\mathrm{k}^{\,4}}{\mathrm{n}}. \tag{21}$$

Changing variables in (20) and integrating over r and $w_\mathrm{i}$, leads to the distribution,

$$f_{n\mathrm{jk}}(\mathrm{n}, \mathrm{x}_\mathrm{i}', \mathrm{x}_\mathrm{j}') = \tag{22}$$
$$\int_{w_\mathrm{i}, \zeta_{\mathrm{ki}}} \frac{\mathrm{n}}{\zeta_\mathrm{k}^{\,4}} f_{\mathrm{x}_\mathrm{i}}(w_\mathrm{i}) f_{\mathrm{x}_\mathrm{k}}(\sqrt{\mathrm{n}^2 - (\mathrm{x}_\mathrm{i}' - w_\mathrm{i})^2 - \mathrm{x}_\mathrm{j}'^2}) f_{\tilde{\mathrm{x}}_\mathrm{i}}(\frac{(\mathrm{x}_\mathrm{i}' - w_\mathrm{i})}{\zeta_\mathrm{k}}) f_{\tilde{\mathrm{x}}_\mathrm{j}}(\frac{\mathrm{x}_\mathrm{j}'}{\zeta_\mathrm{k}}) f_{\tilde{\mathrm{x}}_\mathrm{k}}(\zeta_\mathrm{k}\sqrt{\mathrm{n}^2 - (\mathrm{x}_\mathrm{i}' - w_\mathrm{i})^2 - \mathrm{x}_\mathrm{j}'^2}) dw_\mathrm{i} d\zeta_\mathrm{k}.$$

The symmetry of (19) similarly applies to (22), so that nine equations formed from (14) and (22) form the PDFs of all cases of entering/exiting the box.

## 2.3 Case II: Uniform entrance

### *2.3.1. Opposing faces*

The PDF for the trajectory length when exiting on this facing side can be determined from four independent variables. A line traveling between faces $\mathrm{x}_\mathrm{j} = 0$ and $\mathrm{x}_\mathrm{j}' = \mathrm{X}_\mathrm{j}$ is described by

$$\mathrm{n} = \sqrt{(\mathrm{x}_\mathrm{i}' - \mathrm{x}_\mathrm{i})^2 + \mathrm{X}_\mathrm{j}^{\,2} + (\mathrm{x}_\mathrm{k}' - \mathrm{x}_\mathrm{k})^2}, \tag{23}$$

a problem that resembles the case of lines formed from two uniformly distributed random points on the box surface [12], with the exception that $\mathrm{x}_\mathrm{i}'$ and $\mathrm{x}_\mathrm{k}'$ are no longer uniformly distributed, leaving only two independent random variables. Equation (23) is reduced by substituting $s_\mathrm{i} = (\mathrm{x}_\mathrm{i}' - \mathrm{x}_\mathrm{i})^2$, and

$s_\mathrm{k} = (\mathrm{x}_\mathrm{k}' - \mathrm{x}_\mathrm{k})^2$

$$\mathrm{n} = \sqrt{s_\mathrm{i} + \mathrm{X}_\mathrm{j}^{\,2} + s_\mathrm{k}}, \tag{24}$$

then by (A.4) and (A.6)

$$f_s(s_\mathrm{i}) = \frac{f_{\mathrm{x}_\mathrm{i}}(\mathrm{x}_\mathrm{i}' - \sqrt{s_\mathrm{i}})}{2\sqrt{s_\mathrm{i}}} + \frac{f_{\mathrm{x}_\mathrm{i}}(\mathrm{x}_\mathrm{i}' + \sqrt{s_\mathrm{i}})}{2\sqrt{s_\mathrm{i}}} \tag{25}$$
$$f_s(s_\mathrm{k}) = \frac{f_{\mathrm{x}_\mathrm{k}}(\mathrm{x}_\mathrm{k}' - \sqrt{s_\mathrm{k}})}{2\sqrt{s_\mathrm{k}}} + \frac{f_z(\mathrm{x}_\mathrm{k}' + \sqrt{s_\mathrm{k}})}{2\sqrt{s_\mathrm{k}}}$$

which by (A.1a) and with $S = s_i + s_k$, further reduces to

$$f_S(S+X_j^2) = f_{s_k}(s) \otimes f_{s_i}(s) \cdot \tag{26}$$

Applying (A.5) and using $n = \sqrt{S}$,

$$f_n(n \mid x_i', x_k') = f_S(n^2 + X_j^2)n/2 \tag{27}$$

Applying all even permutations of (i,j,k), provides the PDFs for lines entering through a box side coplanar with one of the Cartesian planes and exiting the opposing face. As in the previous cases, by symmetry

$$f_n(n, x_i', x_k')\,|_{x_j=X_j}^{x_j'=0} = f_n(n, x_i', x_k')\,|_{x_j=0}^{x_j'=X_j}. \tag{28}$$

*2.3.2. Adjacent faces*

Between adjacent sides, the length n is given by

$$n=\sqrt{(x_i'-x_i)^2+x_j'^2+x_k^2}\,. \tag{29}$$

The PDFs can be obtained in the same fashion as the opposing case. Letting $s_i+s_k=(x_i'-x_i)^2+x_k^2$,

by (A.4) the distributions of the independent terms $S_i$ and $S_k$ are,

$$\begin{aligned} f_{s_i}(s_i) &= \frac{f_{x_i}(x_i'-\sqrt{s_i})}{2\sqrt{s_i}} + \frac{f_{x_i}(x_i'+\sqrt{s_i})}{2\sqrt{s_i}} \\ f_{s_k}(s_k) &= \frac{f_{x_k}(\sqrt{s_k})}{2\sqrt{s_k}} + \frac{f_{x_k}(\sqrt{s_k})}{2\sqrt{s_k}}. \end{aligned} \tag{30}$$

then by (A.1a) and using $S=s_i+s_k$,

$$f_S(S) = f_{s_i}(s) \otimes f_{s_k}(s), \tag{31}$$

Applying (A.5) to (31) and shifting by $x_j'^2$

$$f_n(n \mid x_i', x_j') = f_S(n^2 + x_j^2)n/2 \cdot \tag{32}$$

As in the previous section, symmetry there at most six unique PDFs for the adjacent faces,

$$\begin{aligned} f_n(n \mid x_i', x_j')\Big|_{x_j=0}^{x_k'=0} &= f_n(n \mid x_i', x_j')\,|_{x_j=0}^{x_k'=X_k} \\ &= f_n(n \mid x_i', x_j')\,|_{x_j=X_j}^{x_k'=0} \\ &= f_n(n \mid x_i', x_j')\,|_{x_j=X_j}^{x_k'=X_k} \end{aligned} \tag{33}$$

taken over *all* permutations of (i,j,k).

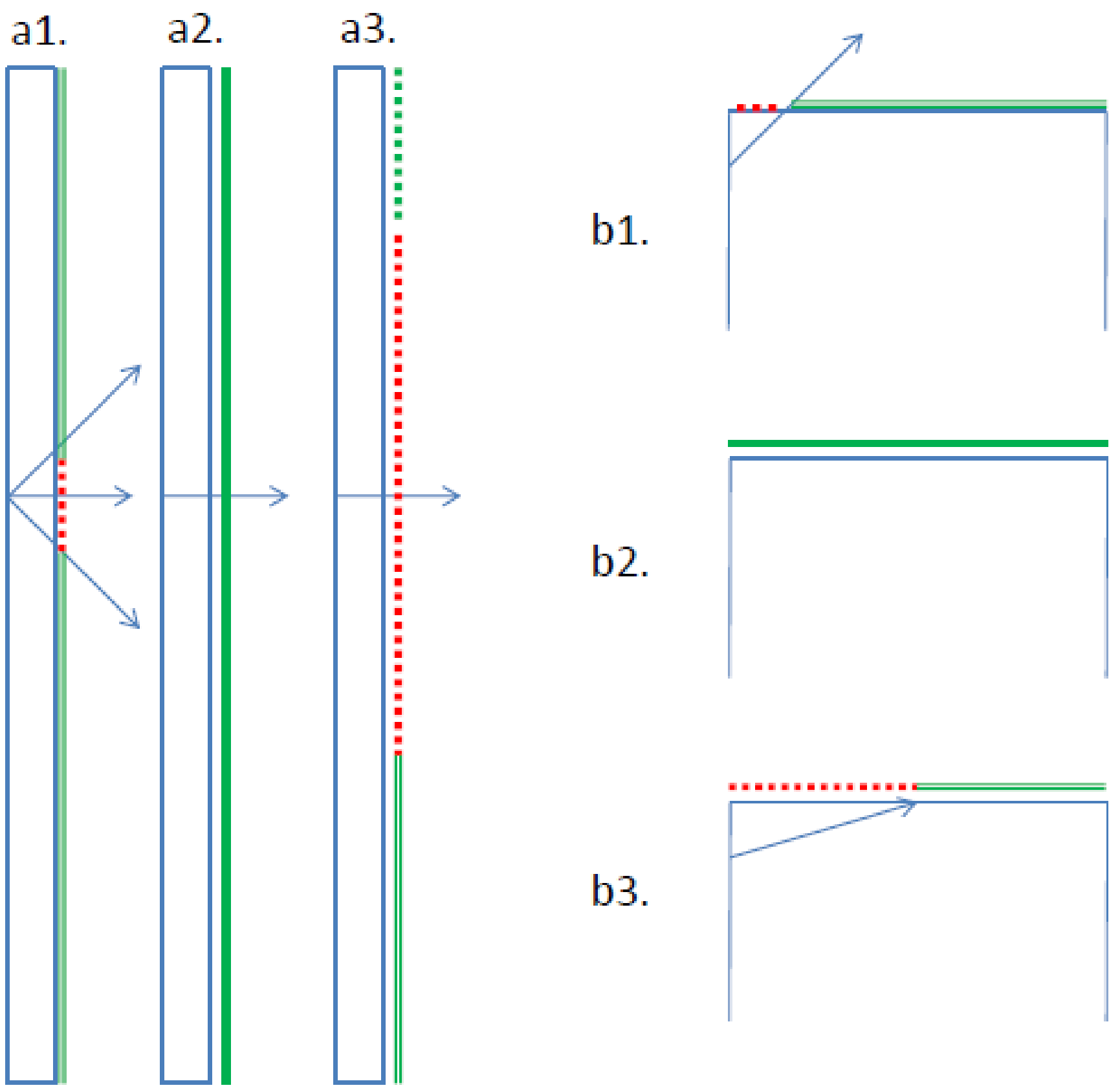

Figure 2. (a) Two-Dimensional illustration of a trajectory across a thin box. From the entry point shown, there is a 50% probability of exiting through the dotted region when (a1) randomness is defined by random orientation and translation (Case I) compared with that of (a2) uniformly distributed random entrance and exit points. The behaviors of both trajectories are contrasted with that of (a3) radial spreading in all directions (Case II). Trajectories exiting an adjacent side are illustrated in (b).

## 3. Examples

### 3.1 Surface Distributions

By way of example, several characteristics of the two cases become apparent. Under the distributions of Case I, a single side of a box might represent a surface where spontaneous emission of a particle is expected with uniform random probability. Such a particle might represent such diverse events as a photon emission, or a single particle from an explosion, whereas Case II might represent acoustic emission emanating from a point on a surface with uniform probability (e.g. from the formation of a bubble) or light intermittently emitted from points on a surface (e.g. an ensemble of distributed fireflies).

Definitions dictating Case I allow for distinct behavior as a function of Box dimension as well as in comparison to Case II. This is in part due to the uniform distribution of trajectory angles, as illustrated in Figure 2. To demonstrate such dependence, distributions were formed by calculating (11), (14), (19), and (22) for three box sizes: a “Short” box (1, 0.1, 1), a “Cubic” box (1, 1, 1) , and a “Long” box (1, 2, 1). To perform these calculations, scripts were written and implemented in Matlab R2010a using Windows 7 64-bit operating system using two quad-core E5-1620 3.6 GHz Xeon processors and 32 GB of RAM, with Matlab typically utilizing 26 GB – 30 GB with running times of 1-2 days. For comparison, 50 million random trajectories were generated for Case I and II, based on calculations of uniform random distributions of relevant variables.

Studied over an entire surface, the distributions can give insight on points to measure, or avoid, in order to achieve maximum/minimum likelihood of encountering a particular length. Such information is particularly relevant when trajectory length is a function of other quantities (radioactive decay, the temperature or velocity of an ejected particle, the amplitude of a spreading wave). An example of the PDF over the opposite and one adjacent face of a cube under Case II is shown with its corresponding randomly-generated distribution in Figure 3.

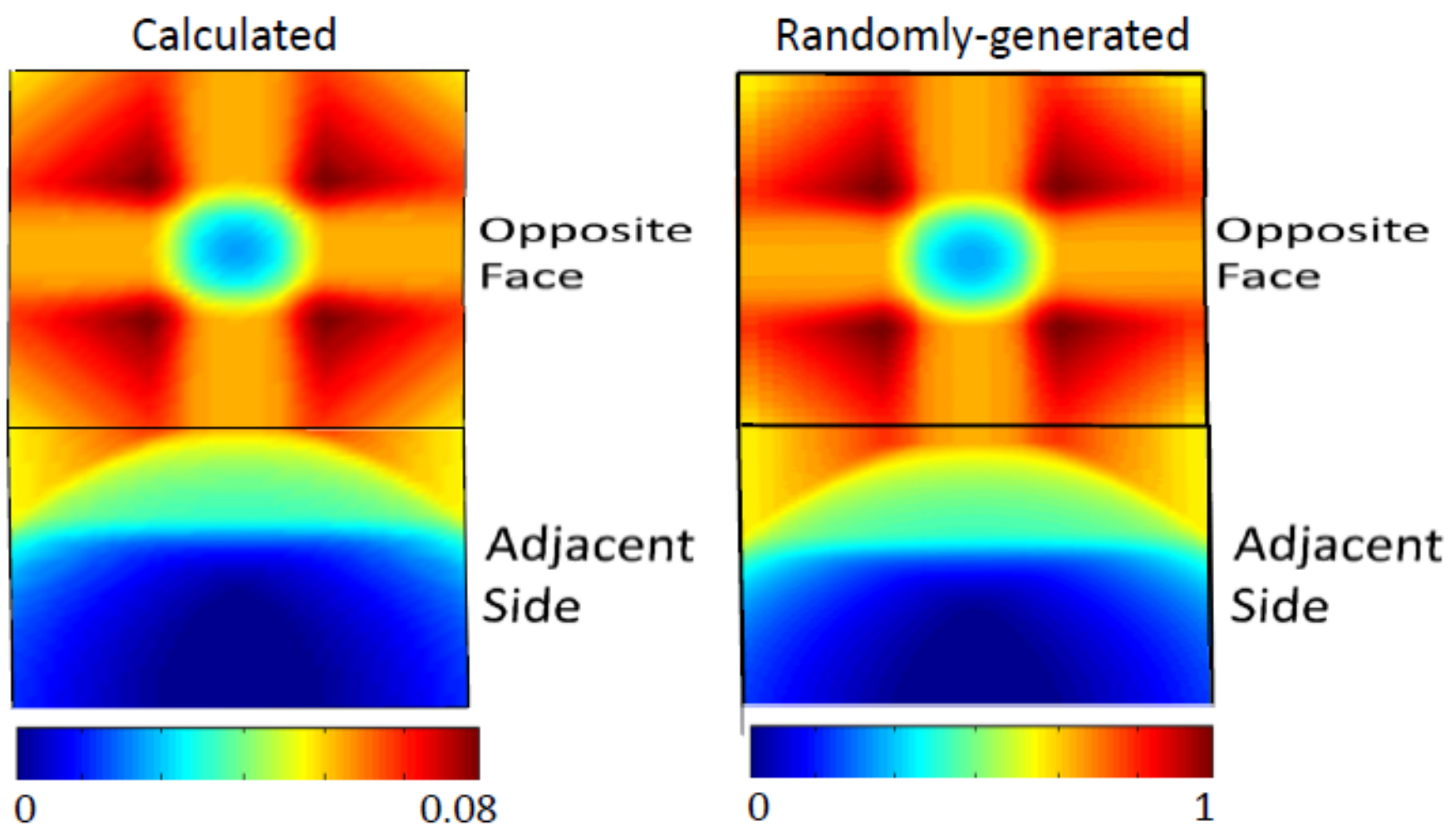

Figure 3. (left) PDF for a trajectory of length 1.17 <= n <= 1.22 through a unit cube and (right) number of trajectories meeting this criteria from a random distribution of 50 million trajectories.

As illustrated in Figure 2b, uniformly distributed rays (Case I) that enter close to a given adjacent side are more likely to exit at lower elevations of that side. This increases the likelihood of measuring trajectories of small length while simultaneously producing the highest concentration of trajectories. Figure 4 displays the randomly generated distributions and PDF indicating the likelihood of a given length on an adjacent side as a function of elevation. The PDF is formed by integrating over one surface variable (here, $x_1'$ ) that, when compared to Case II, reveals a pronounced difference between the two. Case II (Figure 5) contrasts with that of Case I in that the trajectories exiting through low elevations are more likely to be of greater length. Moreover, particle concentration for this case peaks at a higher elevation and with longer length as compared to Case I.

The joint distribution examined over a finite surface element at a given surface location provides information on the trajectory PDF characteristic of the particular location and box dimensions. This might correspond to events recorded by a detector at this location (e.g. a hole of finite size through which a particle may be ejected, or a transducer that records a wave upon exit). Figure 6 shows an example of the PDF for Case II at a given location $(x_1' = 1; x_3' = 0.25)|_{x_2' = X_2}$ over a small area ($\Delta_1' = \Delta_3' = 0.003$), along with a corresponding randomly-generated distribution. Such a distribution provides a signature behavior whose characteristics are particular to the specific location, type of randomness, and box dimensions.

### 3.2 Combined distribution

Distributions of trajectory length over an entire box can be expressed by the combined density functions of all sides. For Case I, the distribution is obtained by integrating (14) over $X_i'$ and $X_k'$ and (22) over $X_i'$ and $X_j'$. Then by (10), (11), (18), and (19), the number of terms is reduced to 9 from a possible 30 that may be summed to form the PDF

$$f_{\mathrm{I}}(\mathrm{n}) = \sum_{3} 2 f_{\mathrm{n}_{\parallel}}(\mathrm{n}) \Big|_{x_j=0}^{x'_j=X_j} \mathrm{P}_{x_j=0} + \sum_{6} 4 f_{\mathrm{n}_{\perp}}(\mathrm{n}) \Big|_{x_j=0}^{x'_j=X_j} \mathrm{P}_{x_j=0}, \tag{34}$$

with $f_{\mathrm{n}_{\parallel}}$ denoting PDFs of trajectories through opposite faces and $f_{\mathrm{n}_{\parallel}}$ through adjacent sides. Likewise, from (28) and (33)

$$f_{\mathrm{II}}(\mathrm{n}) = \sum_{3} 2 f_{\mathrm{n}_{\parallel}}(\mathrm{n}) \Big|_{x_j=0}^{x'_j=X_j} \mathrm{P}_{x_j=0} + \sum_{6} 4 f_{\mathrm{n}_{\perp}}(\mathrm{n}) \Big|_{x_j=0}^{x'_j=X_j} \mathrm{P}_{x_j=0}. \tag{35}$$

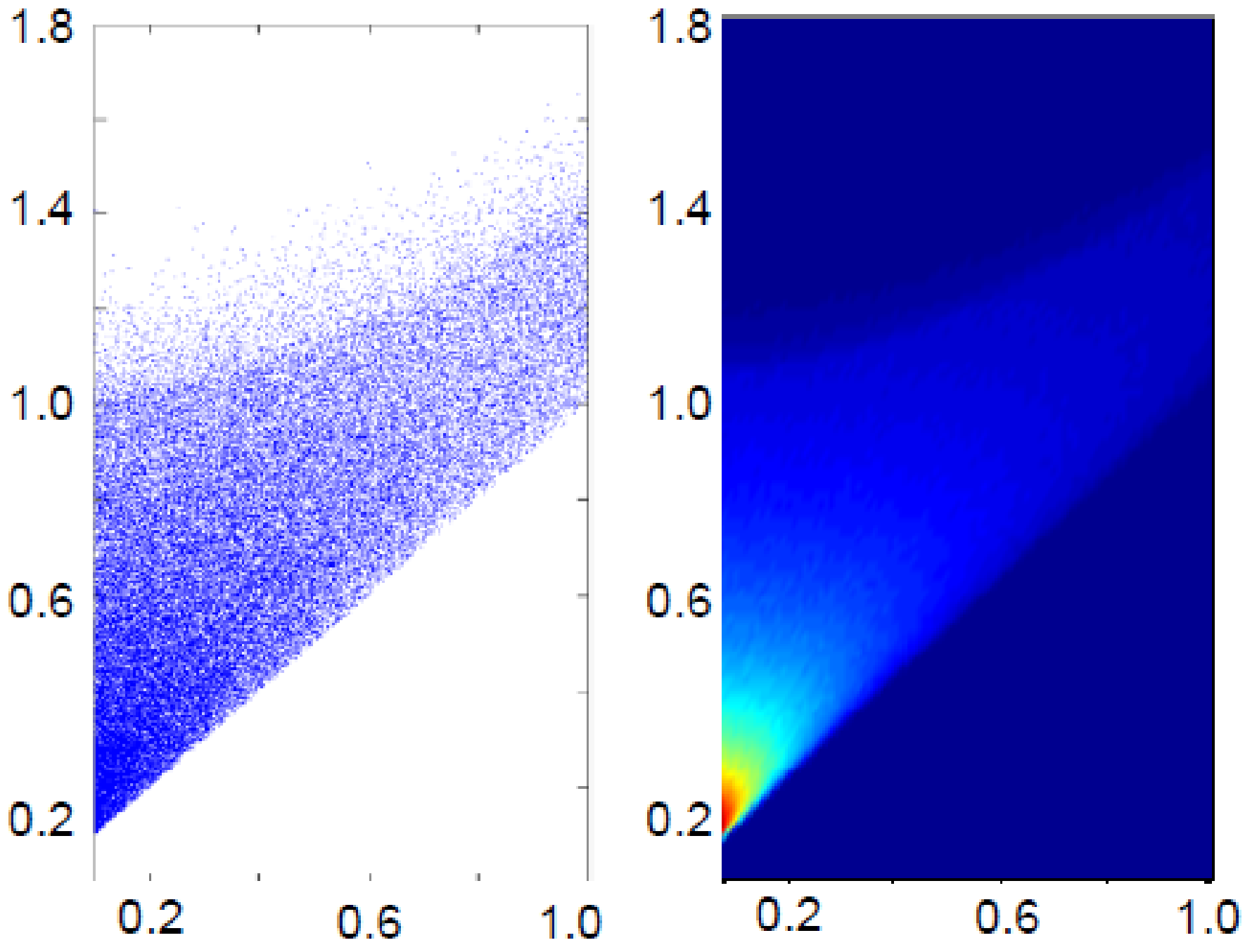

Figure 4. Distribution of trajectory lengths as a function of exit position on an adjacent axis of a cube for Case I showing (left) 31,000 of 50 million trajectories, and its associated two-dimensional distribution of the same data (right).

Resulting equations (34) and (35) may be used to find expectation values for trajectory length through a box, which has been both described and applied in past works [2], [5], [9], [12]. For physical problems examining emissions from a single surface, or equivalently particles entering from a single opening, the relevant combined distribution reduces to the sum of the surface's four adjacent sides and opposite side. Under this scenario PDFs and their associated randomly-generated distributions are shown in Figure 7 for the Short case, Figure 8 for the Cubic case, and figure 9 for the Long case.

## 4. Concluding Remarks

Distributions of the distance between two random points in and on a rectangular cuboid (box) surface has been well studied in area of geometric probability [9]. Varying methods have been used to describe the PDFs and expectation values between such points, which depend upon the particular type of randomness being studied [4]. These values form geometric constants that are of potential practical interest in acoustics, particle and atomic physics, radiology, material science, and a number of other areas. Such studies have been expanded to the more general case of a rectangular box of arbitrary length, width and height [5].

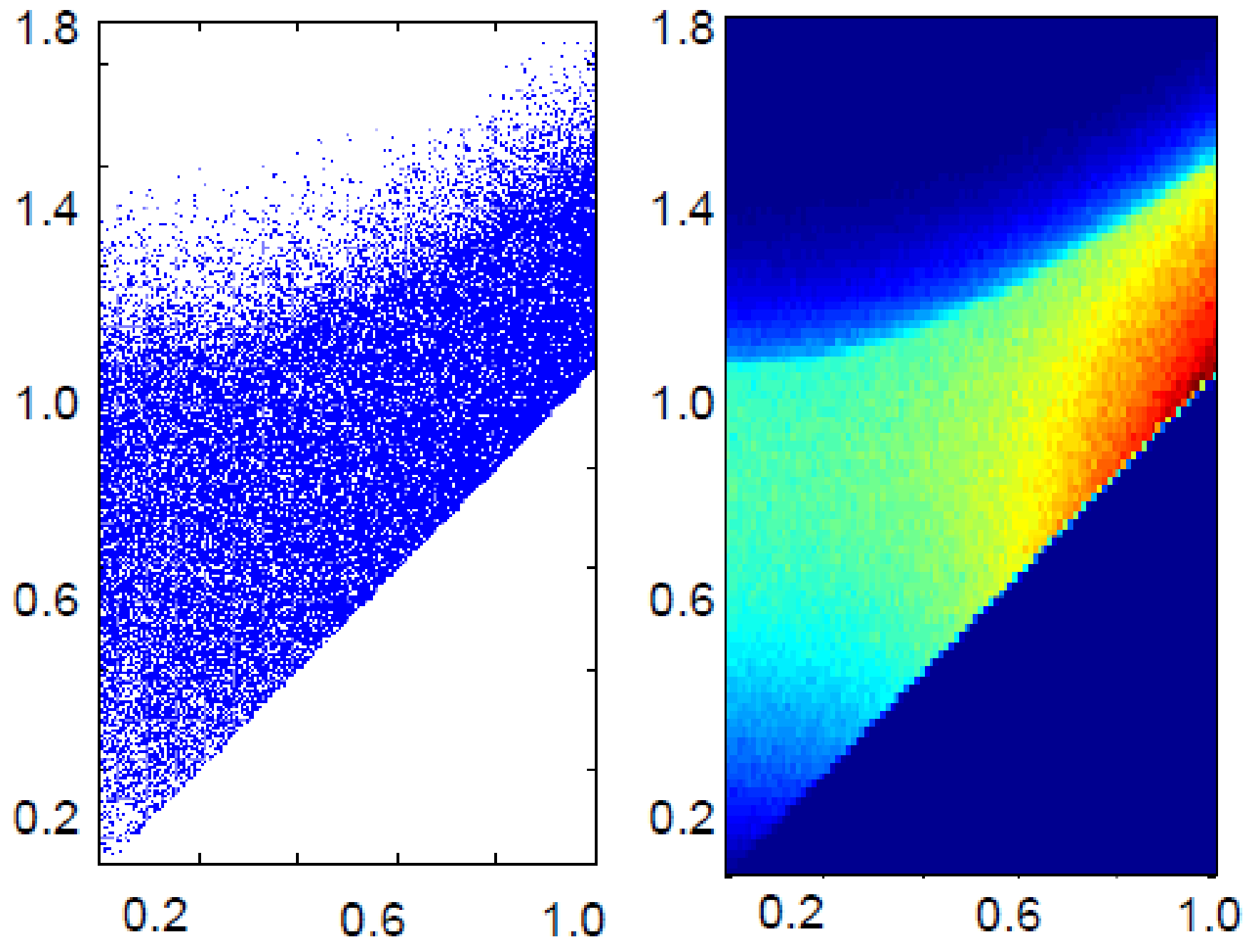

Figure 5. Distribution of trajectory lengths as a function of exit position on an adjacent axis of a cube for Case II showing (left) 31,000 of 50 million trajectories, and its associated two-dimensional distribution (right).

For the purposes of the present study, two random points on a box surface were considered as trajectory intersections with the box. The joint probability density was calculated for the length and exit position of two types of randomness, including rays traveling in a uniformly-random direction through a box of arbitrary relative dimension and a radial vector extending from a random point on the box surface. These cases may represent a number of physical situations, perhaps the most interesting being interpretation as a particle or wave, respectively. Our own ongoing efforts will continue toward applying such analysis toward signal evaluation in ultrasound.

## Acknowledgements

This work was supported by the Cleveland Clinic Foundation. A portion of the work was motivated by and funded under grant R01EB014296 from the NIH National Institute of Biomedical Imaging and Bioengineering. The content is solely the responsibility of the author.

## References

1. E. W. WEISSTEIN, Geometric Constants *Wolfram MathWorld*." Available: http://mathworld.wolfram.com/topics/GeometricConstants.html.
2. D. SHI, A probabilistic approach to estimate particle size distribution without shape assumption. Part I: Theory and Testing on spherical particles, *Materials Characterization*, 27(1): 35–44 (1991).
3. J.-H. HAN and D.-Y. KIM, Determination of three-dimensional grain size distribution by linear intercept measurement, *Acta Materialia*, 46(6):2021–2028 (1998).
4. R. COLEMAN, Random Paths through Convex Bodies, *Journal of Applied Probability*, 6(2):430–441 (1969).
5. R. COLEMAN, Intercept Lengths of Random Probes through Boxes, *Journal of Applied Probability*, 18(1):276–282 (1981).

6. H. J. PALTIEL, C. R. ESTRADA, A. I. ALOMARI, C. STAMOULS, C. C. PASSEROTTE, F. C. MERAL, R. S. LEE, and G. T. CLEMENT, Multiplanar Dynamic Contrast-enhanced US Assessment of Blood Flow in a Rabbit Model of Testicular Torsion, *Ultrasound in Medicine & Biology*, IN PRESS, (2013).
7. H. J. PALTIEL, H. M. PADUA, P. C. GARGOLLO, G. M. CANNON, A. I. ALOMARI, R. YU, and G. T. CLEMENT, Contrast-enhanced, real-time volumetric ultrasound imaging of tissue perfusion: preliminary results in a rabbit model of testicular torsion, *Physics in Medicine and Biology*, 56:2183-2197 (2011).
8. F. C. MERAL and G. T. CLEMENT, 128 Element ultrasound array for transcranial imaging, *2010 IEEE Ultrasonics Symposium (IUS)* pp. 1984–1987 (2010).
9. A. M. MATHAI, P. MOSCHOPOULOS, and G. PEDERZOLI, Distance between random points in a cube, *Statistica*, 59(1): 61-81 (1999).
10. J. PHILIP, The probability distribution of the distance between two random points in a box, *TRITA MAT* 7(10) 2007.
11. G. MARSAGLIA, Choosing a Point from the Surface of a Sphere, *The Annals of Mathematical Statistics*, 43(2):645–646 (1972).
12. J. PHILIP, Calculation of Expected Distance on a Unit Cube, *www.math.kth.se/~johanph,* 2007.

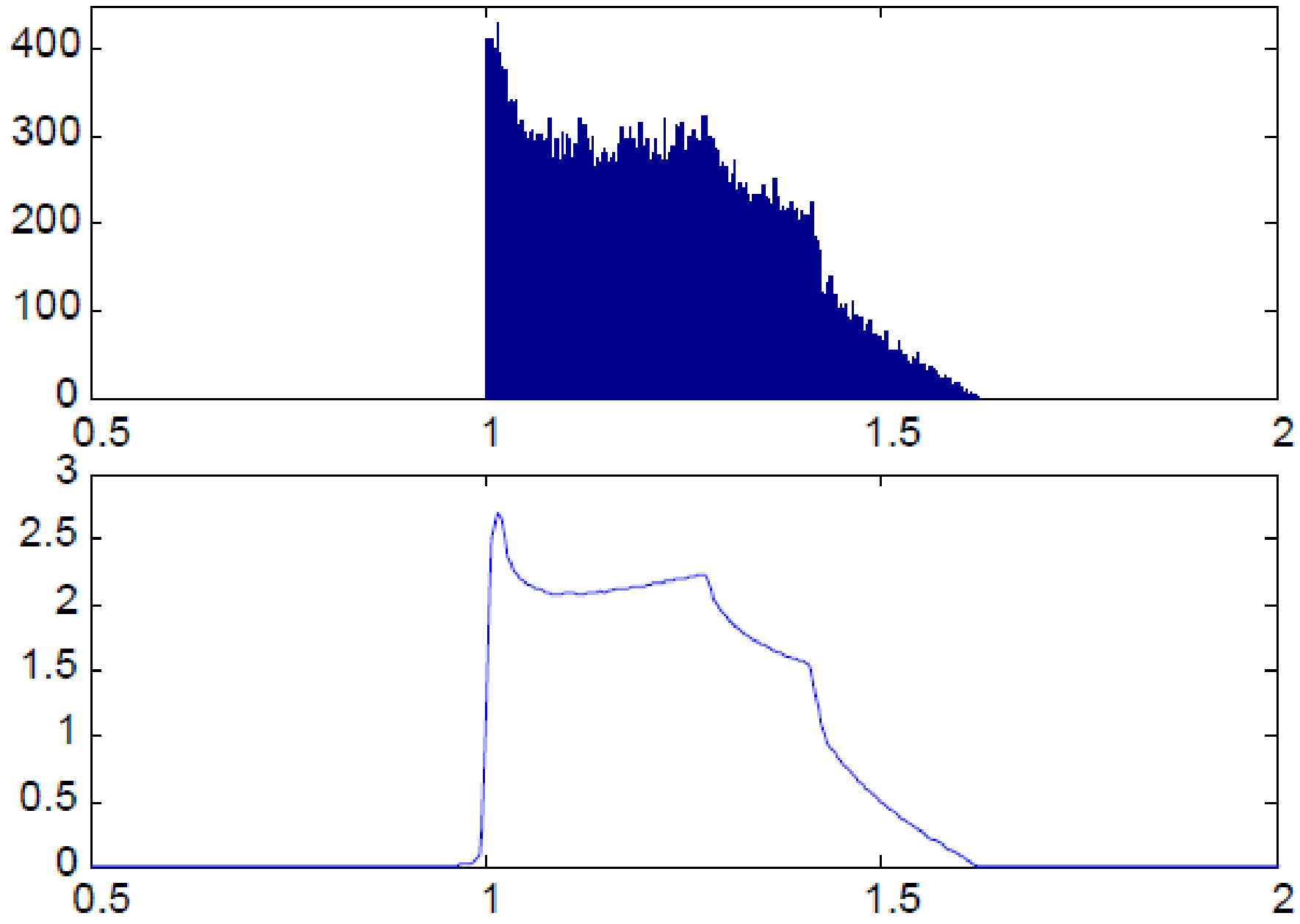

Figure 6. (top) Histogram of trajectory lengths of randomly-generated trajectories following Case II examined at surface location (xOUT=1, yOUT=1, zOUT=0.25 ) on a unit cube and (bottom) the calculated PDF at the location.

## Appendix A. Some useful PDF transformations

Let $x$ and $y$ be uniformly distributed independent random variables. Distributions formed from the sums or products of these variables relevant to the present study are summarized here.

$s = x + y$:

$$f(s) = f_x(s) \otimes f_y(s) \tag{A.1a}$$

$$s = x - y: \qquad \text{(A.1b)}$$
$$f(\mathrm{s}) = f_x(\mathrm{s}) \otimes f_y(-\mathrm{s})$$

$$s = x / y:$$

Let w=y. The Jacobian is $\begin{vmatrix} \hat{\mathbf{i}} & \hat{\mathbf{j}} & \hat{\mathbf{k}} \\ \frac{\partial s}{\partial x} & \frac{\partial s}{\partial y} & 0 \\ \frac{\partial w}{\partial x} & \frac{\partial w}{\partial y} & 0 \end{vmatrix} = \frac{1}{y}\hat{\mathbf{k}}$ (A.2)

$$f(\mathrm{s,w}) = f_x(\mathrm{s\,w})\, f_y(\mathrm{w})$$
$$f(\mathrm{s}) = \int_0^\infty f_x(\mathrm{s\,w})\, f_y(\mathrm{w})\, \mathrm{w\, dw}$$

$$s = xy:$$

Let w=y. The Jacobian is $\begin{vmatrix} \hat{\mathbf{i}} & \hat{\mathbf{j}} & \hat{\mathbf{k}} \\ \frac{\partial s}{\partial x} & \frac{\partial s}{\partial y} & 0 \\ \frac{\partial w}{\partial x} & \frac{\partial w}{\partial y} & 0 \end{vmatrix} = y\hat{\mathbf{k}}$ (A.3)

$$f(\mathrm{s,w}) = \frac{f_x(\mathrm{s/w})\, f_y(\mathrm{w})}{\mathrm{w}}$$
$$f(\mathrm{s}) = \int_0^\infty \frac{f_x(\mathrm{s/w})\, f_y(\mathrm{w})}{\mathrm{w}} \mathrm{dw}$$

$$s = x^2: \qquad \text{(A.4)}$$
$f(\mathrm{s}) = \frac{f_x(\sqrt{\mathrm{s}})}{2\sqrt{s}} + \frac{f_x(-\sqrt{\mathrm{s}})}{2\sqrt{s}}$ i.e. Neg part of $f_x$ contributes to pos $f$

$$s = \sqrt{x}: \qquad \text{(A.5)}$$
$$f(\mathrm{s}) = f_x(\mathrm{s}^2)\,\mathrm{s}/2$$

## Appendix B. Jacobians for the joint distributions of Case I and Case II

Case 1:

$$
J=\begin{vmatrix}
1 & 0 & 0 & 0 & 0 \\
0 & 1 & 0 & -X_j\frac{\tilde{x}_k}{\tilde{x}_j^{\,2}} & -\frac{X_j}{\tilde{x}_j} \\
1 & 0 & X_j\frac{1}{\tilde{x}_j} & -X_j\frac{\tilde{x}_i}{\tilde{x}_j^{\,2}} & 0 \\
0 & 0 & X_j\frac{\tilde{x}_i}{\tilde{x}_j\sqrt{\tilde{x}_i+\tilde{x}_j+\tilde{x}_k}} & \frac{X_j}{\sqrt{\tilde{x}_i+\tilde{x}_j+\tilde{x}_k}}-X_j\frac{\sqrt{\tilde{x}_i+\tilde{x}_j+\tilde{x}_k}}{\tilde{x}_y^{\,2}} & X_j\frac{\tilde{x}_k}{\tilde{x}_j\sqrt{\tilde{x}_i+\tilde{x}_j+\tilde{x}_k}} \\
0 & 0 & \frac{\tilde{x}_i}{\sqrt{\tilde{x}_i+\tilde{x}_j+\tilde{x}_k}} & \frac{\tilde{x}_j}{\sqrt{\tilde{x}_i+\tilde{x}_j+\tilde{x}_k}} & \frac{\tilde{x}_k}{\sqrt{\tilde{x}_i+\tilde{x}_j+\tilde{x}_k}}
\end{vmatrix}
$$

(columns: $x_i$, $x_k$, $\tilde{x}_i$, $\tilde{x}_j$, $\tilde{x}_k$; rows: $w_i$, $x_k'$, $x_i'$, n, r)

(B.1)

$$
=-\frac{X_j^{\,2}x_k}{\tilde{x}_j^{\,3}}=\frac{n^2\sqrt{n^2-(x_i'-w_i)^2-X_j^{\,2}}}{X_j r^2}
$$

Case 2:

$$
J=\begin{vmatrix}
1 & 0 & 0 & 0 & 0 \\
0 & \frac{1}{\tilde{x}_k} & 0 & 0 & -\frac{x_k}{\tilde{x}_k^{\,2}} \\
0 & \frac{\tilde{x}_j}{\tilde{x}_k} & 0 & \frac{x_k}{\tilde{x}_k} & -\frac{x_k\tilde{x}_j}{\tilde{x}_k^{\,2}} \\
1 & \frac{\tilde{x}_i}{\tilde{x}_k} & \frac{x_k}{\tilde{x}_k} & 0 & -\frac{x_k\tilde{x}_i}{\tilde{x}_k^{\,2}} \\
0 & \frac{\sqrt{\tilde{x}_i+\tilde{x}_j+\tilde{x}_k}}{\tilde{x}_k} & \frac{x_k}{\tilde{x}_k\sqrt{\tilde{x}_i+\tilde{x}_j+\tilde{x}_k}} & \frac{x_k\tilde{x}_i}{\tilde{x}_k\sqrt{\tilde{x}_i+\tilde{x}_j+\tilde{x}_k}} & \frac{x_k}{\sqrt{\tilde{x}_i+\tilde{x}_j+\tilde{x}_k}}-\frac{x_k\sqrt{\tilde{x}_i+\tilde{x}_j+\tilde{x}_k}}{\tilde{x}_k^{\,2}}
\end{vmatrix}
$$

(columns: $x_i$, $x_k$, $\tilde{x}_i$, $\tilde{x}_j$, $\tilde{x}_k$; rows: $w_i$, $\zeta_k$, $x_j'$, $x_i'$, n)

(B.2)

$$
=-\frac{x_k^{\,3}}{\sqrt{\tilde{x}_i+\tilde{x}_j+\tilde{x}_k}}=\frac{\zeta_k^{\,3}}{r}
$$

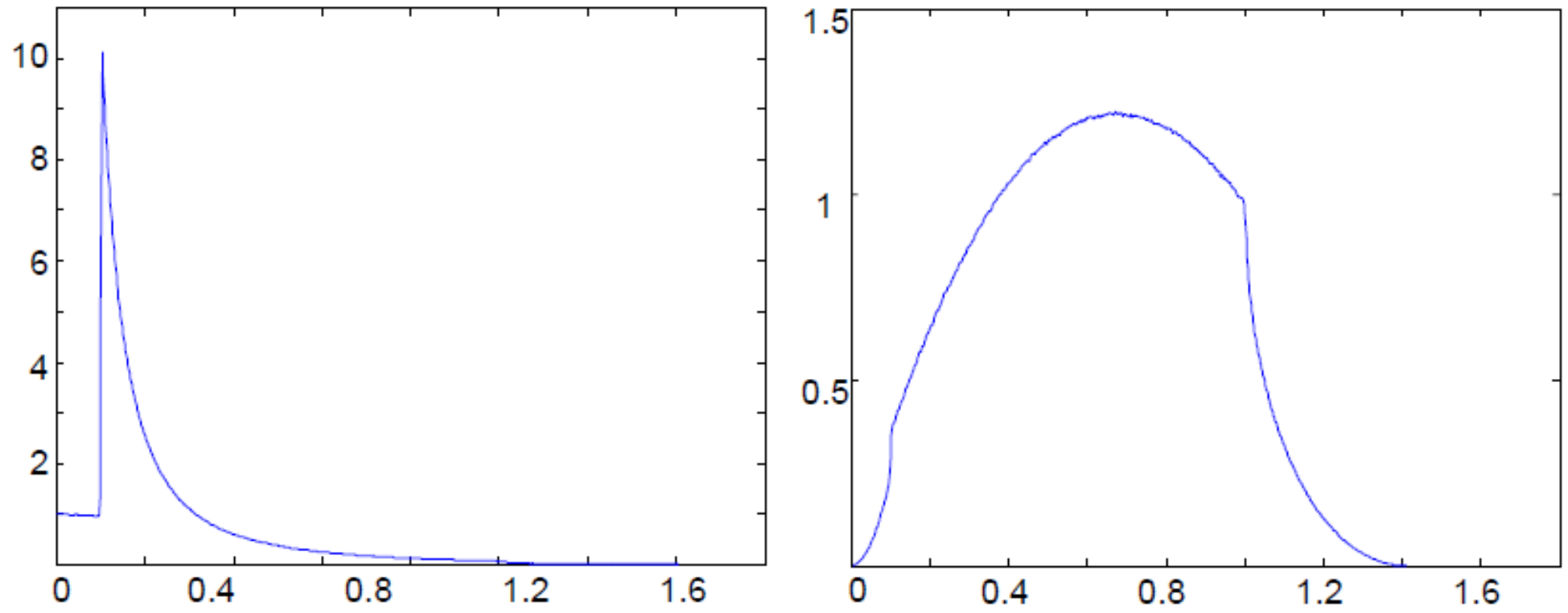

Figure 7. Combined PDFs for Case I (left) and Case II (right) with box dimensions (1,0.1,1) and trajectories entering through a single (1,1) side.

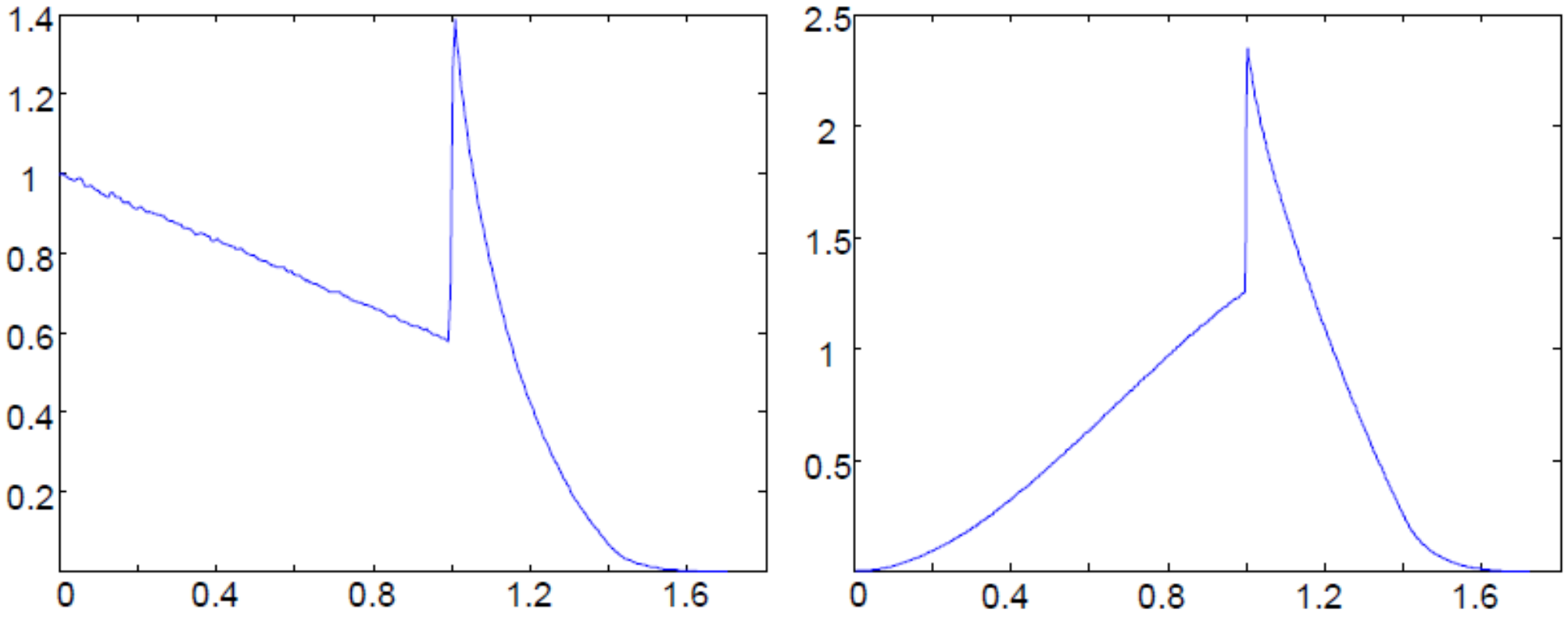

Figure 8. Combined PDFs for Case I (left) and Case II (right) with box dimensions (1,1,1) and trajectories entering through a single (1,1) side.

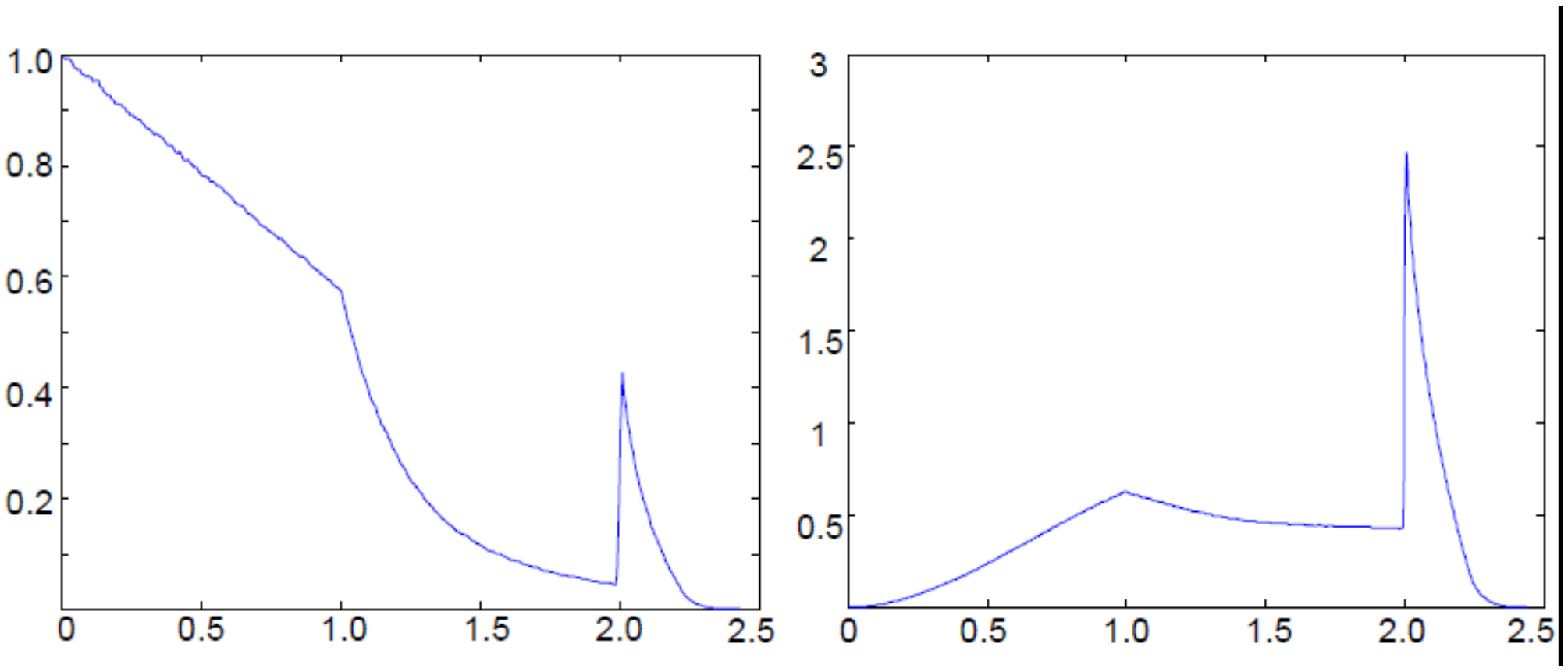

Figure 9. Combined PDFs for Case I (left) and Case II (right) with box dimensions (1,2,1) and trajectories entering through a single (1,1) side.